\documentclass{amsart}
\usepackage{amssymb}

\usepackage{tikz}

 \theoremstyle{definition}
 
 \theoremstyle{remark}

 \numberwithin{equation}{section}

\begin{document}
%
%
%

\title
 {An Inaccessible Graph}
 
\author{M. J. Dunwoody}

\address{%
University of Southampton\\
Southampton\\
SO16 7GR\\
U.K}

\email {M.J.Dunwoody@soton.ac.uk}



\begin{abstract}
An inaccessible, vertex transitive, locally finite graph is described.   This graph is not quasi-isometric to a Cayley graph. 
\end{abstract}

\maketitle

\section{Introduction}

\newcommand{\R}{\mathbb{R}}

\newcommand{\Z}{\mathbb{Z}}
\newcommand{\cd}{\mathcal{D}}
\newcommand{\cf}{\mathcal{F}}

\newcommand{\cn}{\mathcal{N}}
\newcommand{\cg}{\mathcal{G}}

\newcommand{\cc}{\mathcal{C}}
\newcommand{\cm}{\mathcal{M}}
\newcommand{\ct}{\mathcal{T}}
\newcommand{\cs}{\mathcal{S}}
\newcommand{\ce}{\mathcal{E}}
\newcommand{\cb}{\mathcal{B}}
\newcommand{\cp}{\mathcal{P}}
\newcommand{\cv}{\mathcal{V}}
\newcommand{\cl}{\mathcal{L}}
\newcommand{\cw}{\mathcal{W}}
\newcounter{fig}
\setcounter{fig}{0}

Let $X$ be a locally finite connected graph.
A {\it ray} is a sequence of distinct vertices $v_0, v_1, \dots  $ such that $v_i$ is adjacent to $v_{i+1}$ for each $i = 1, 2, \dots $.  Obviously
for a ray to exist, the graph $X$ has to be infinite.
For any two vertices $u, v \in VX$ let $d(u,v)$ be the length of a shortest path joining $u, v$.

We say that two rays $R, R'$  belong to the same {\it end} $\omega $, if for no finite subset $F$ of $VX$ or $EX$  do $R_1$
and $R_2$ eventually lie in distinct components of $X \setminus F$.
We define $\ce (X)$ to be the set of ends of $X$.

We say that $\omega $ is {\it thin } if it does not contain infinitely many vertex disjoint rays
   As in \cite{[TW]} the end $\omega $ is said to be {\it thick} if it is not thin.

In their nice paper \cite{[TW]} Thomassen and Woess define an accessible graph.  
A graph $X$ is {\it accessible} if there is some natural number $k$ such that for any two ends $\omega _1$ and $\omega _2$ of $X$, there is a set 
$F$ of at most $k$ vertices in $X$ such that $F$ separates $\omega _1$ and $\omega _2$, i.e. removing $F$ from $X$ disconnects the graph
in such a way  that   rays $R_1,  R_2$ of $\omega _1, \omega _2$ respectively eventually lie in distinct components of $X\setminus F$.

A finitely generated group $G$ is said to have more than one end ($e(G) > 1$) if its Cayley graph $X(G, S)$ with respect to a finite generating set $S$
has more than one end.  This property is independent of the generating set $S$ chosen.  Stallings \cite{[S1971]} showed that 
if $e(G) > 1$ then $G$ {\it splits} over a finite subgroup, i.e. either $G = A*_CB$ where $C$ is finite, $C \not= A, C\not= B$ or
$G$ is an HNN extension $G = A*_C = \langle A, t| t^{-1}ct = \theta (c)\rangle$, where $C$ is finite, $C \leq A$ and $\theta : C \rightarrow A$ is an injective homomorphism.
A group is {\it accessible } if the process of successively factorizing factors that split in a decomposition of $G$  eventually 
terminates  with factors that are finite or one ended.

Thomassen and Woess show that  the Cayley graph of a finitely generated group $G$  is accessible if and only if  $G$ is accessible.
In\cite{[D2], [D3]} I have given examples of inaccessible groups, and so not every locally finite connected graph is accessible.

Let $\omega $ be an end of $X$.   As in \cite{[TW]}, p259 define $k(\omega )$ to be the smallest integer $k$ such that $\omega $ can be separated
from any other end by at most $k$ vertices.  If this number does not exist, put $k(\omega ) = \infty$.

Thomassen and Woess show that $X$ is accessible if and only if $k(\omega ) < \infty$ for every end $\omega $.
We say that an end $\omega $ is {\it special} if $k(\omega ) = \infty$.

In this paper we construct a locally finite, connected, inaccessible, vertex transitive  graph $X$.
The property of being inaccessible is  invariant under quasi-isometry.
If $X, Y$ are graphs, then a  quasi-isometry $\theta : X \rightarrow Y$ induces a bijection $\ce (\theta ) : \ce (X) \rightarrow \ce (Y)$
which takes thick ends to thick ends, and special ends to special ends. One can put a topology on $\ce (X)$ in a natural way.
The map $\ce (\theta)$ is then a homeomorphism.

Woess asked in \cite{[W1], [SW]} if every vertex transitive, locally finite graph is quasi-isometric to a Cayley graph. 
It was shown in \cite{[EFW1], [EFW2]} that the Diestel-Leader graph  $DL(m,n), m\not= n$  (see \cite{[DL]} or \cite{[W1]}) is not quasi-isometric to a Cayley graph, answering the question of Woess.
It is shown here that the graph $X$  is another example.
I originally thought that $X$ was hyperbolic, and the fact that $X$ was not quasi-isometric to a Cayley graph then followed because
a hyperbolic group is finitely presented, and would therefore have an accessible Cayley graph by \cite{[D1]}.  However  there are arbitrarily large cycles in $X$ for which the distance apart of two vertices in the cycle is the same as that in $X$.
This cannot happen in a hyperbolic graph.  It seems likely that a hyperbolic graph must be accessible.

The vertex transitive graph $X$ we construct  is based on a construction  in \cite{[DJ1]}.  In that paper, Mary Jones and I construct a finitely generated group $G$ for which $G \cong A * _C G$ where $C$ is infinite cyclic.
The vertex set of the graph $X$ is the set  of left  cosets of $D$ in $G$, where $D$ has index $2$ in $C$.    One could take the vertex set of
$X$ to be the left cosets of $A$ or $C$ as they are commensurable with $D$.   In fact it is easier to work with a $G$-graph $Y$ quasi-isometric
to $X$, in which there are two orbits of vertices for the action of $G$ on $Y$.

 In general, if a group  $G$ is the commensurizer
of a subgroup $H$,  and $G$
is generated by $H\cup S$, then one can construct a  vertex transitive, connected graph, in which the vertices are the cosets of $H$, and there are edges $(H, sH)$ for each $s \in S$.  If
$G$ actually normalizes $H$, then this graph is a Cayley graph for $G/H$.  Conversely if  $X$ is a connected, vertex transitive, locally finite graph and $H$ is
the stabilizer of a vertex $v$, then $G$ is the commensurizer of $H$ and $G$ is generated by $H\cup S$, where $S$ is any subset of $G$ with the property that
 for each $u$ adjacent to $v$ there is an $s\in S$ such that $sv = u$.

 The graph $Y$ has an orbit of cut points, i.e. vertices whose removal disconnects the graph.   It is well known that cut points in a graph
 give rise to a tree decomposition.      This is described - for example -   in \cite{[DK]},  in which the theory of structure trees is extended to 
 graphs that can be disconnected by removing finitely many vertices rather than finitely many edges.
 The cut point tree $T$ for $Y$  has two orbits of vertices under $G$.    One orbit corresponds to the set of  $2$-blocks, where each $2$-block is 
a maximal $2$-connected subgraph, and the other orbit corresponds to the cut points.
It is then shown that after  a subdivision and two folding operations, each of which is a quasi-isometry, and removing spikes (a spike is an edge with a vertex of degree one) each $2$-block becomes a graph isomorphic
to $Y$.   Thus the graph $Y$ has a self-similarity property that comes from the fact that $G \cong A * _C G$ where $C$ is infinite cyclic.
One would not expect this to happen in a Cayley graph, as it is not possible that for a finitely generated group $G$ to be isomorphic
to $A*_CG$ where $C$ is finite.   This follows from a result of Linnell \cite{[PAL]}, which indicates that in  a process of successively factorizing factors that split in a decomposition of an inaccessilbe group $G$,  the size of the finite groups over which the factors split must increase.

Thus after carrying out the subdivision and folding operations, the graph $Y=Y_1$ becomes a graph  $Y_2$ which has a single orbit of disconnecting edges.  Removing (the interior of) all these
edges will give a single orbit of points each with stabilizer a conjugate of $A$, and a second orbit, consisting  of $2$-blocks each of which is isomorphic 
to $Y$, with stabilizer conjugate to the subgroup of $G$ which is the second factor in the decomposition $G \cong A * _C G$.
If we repeat this process $n-1$ times, then we a obtain 
a graph $Y_n$ which has $n-1$  orbits of disconnecting edges.    Removing these edges produces  $n-1$ orbits of vertices each
of which has finite stabilizer, isomorphic to $A$, and  a single orbit of $2$-blocks   each of which is isomorphic  to $Y$.    Let $B_n$ be one of these blocks.
The graph $Y$ has an orbit of subgraphs each of which is a  trivalent tree.   Let $Z$ be a particular trivalent subtree of $Y$.   
Although the folding operations do involve folding $Z$,  the result of the operations is another trivalent tree.  We will see that any two rays in $Z$ represent a particular special end  $\omega $  of $Y$.
There will also be uncountably many special ends that do not correspond to a translate of $Z$.   A  ray representing a special end must eventually lie in a translate of
$B_n$, since otherwise it will represent a thin end.   However the initial number $x_n$ of points in the ray  outside  a translate of $B_n$ may tend to infinity with $n$.   There
will be uncountably many such special ends.   If the ray eventually ends up in a translate of $Z$, then $x_n$ is bounded, since each
translate of $Z$ lies in a translate of $B_n$.    Since each translate of $B_n$ contains a translate of $Z$, the orbit of $\omega $ is dense in the
space of special ends.  

We will show that in a Cayley graph, if there is a countable set of special ends which is dense in the subspace of all special ends, then there must be
 a special end corresponding to a $1$-ended subgraph.
 There is no special end of $Y$ corresponding to a $1$-ended subgraph, and so the graph $Y $ cannot be quasi-isometric to 
a Cayley graph.



As it is important in our construction, we repeat the description
of $G$ below.
In another paper  \cite{[DJ2]}, Mary Jones and I went on to construct a finitely generated group $G_1$ for which $G_1 \cong G_1 *_{C_1} G_1$ with $C_1$ infinite cyclic.   It might be expected that the coset graph $X_1$
of $C_1 $ in $G_1$ has similar properties to $X$.   This will not be the case.   Although $X_1$ is inaccessible and locally finite, it is quasi-isometric
to a Cayley graph.   This is because $C_1$ contains a central subgroup $Z$ as a subgroup of finite index.  Then $X_1$ is quasi-isometric
to the Cayley graph of $G_1/Z$.   
\section{The graph}

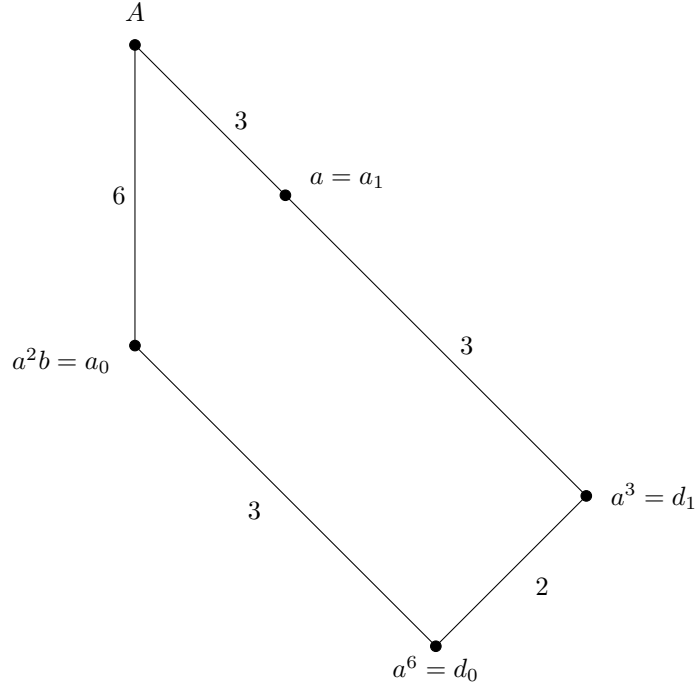
\begin{figure}

\centering
\begin{tikzpicture}[scale=2]

    \path (0,0) coordinate (p1);
    \path (0, 2) coordinate (p2);
    \path (3, -1) coordinate (p3);
    \path (2, -2) coordinate (p4);
    \path (1,1) coordinate (p5);
    \path (1,-1) coordinate (p6);
   \filldraw (p1) circle (1pt) ;
   
     \filldraw (p2) circle (1pt);
  \filldraw (p3) circle (1pt);  
 
  \filldraw (p4) circle (1pt);
    \filldraw (p5) circle (1pt) ;

     \draw (p1) -- (p2) -- (p3) --(p4) -- (p1);
     \draw (.9,-1.1) node [left] {$3$ };
  \draw (2.1,0) node [right] {$3$ };
  \draw (.6, 1.5) node [right] {$3$ };
  \draw (0, 1) node [left] {$6$ };
\draw (2.6, -1.6) node [right] {$2$} ;
\draw (-.1,-.1) node [left] {$a^2b = a_0$} ;
\draw (0, 2.1) node [above] {$A$} ;
\draw (3.1, -1) node [right] {$a^3 = d_1$} ;
\draw (2,-2) node [below] {$a^6 = d_0$} ;
\draw (1.1,1.1) node [right] {$a = a_1$} ;

\end{tikzpicture}
\vskip 1cm \caption{Subgroup lattice in $A$}\label{Fig 1}
\end{figure}

 We recall the group $G$ constructed in \cite{[DJ1]}.  Let $A= \langle a , b | b^3 =1, a^{-1}ba = b^{-1}\rangle$.  As noted in \cite{[DJ1]}, $a^2$ is in the centre of $A$ and $A/\langle a^2 \rangle \cong S_3$.
Also $A$ is generated by $a^3$ and $a^2b$ since $a^{-3}(a^2b)a^3 =a^2b^{-1}$, and so $b^2 = b^{-1} \in \langle a^3, a^2b\rangle$.
The group $A$ has a lattice of subgroups as in Fig \ref{Fig 1}.

Put $x = a^3, y = a^2b$.  Then, since $a^2$ is central $x^2 = y^3$.   Also $y^{-1}x = y^2x^{-1} = b^{-1}a$ and $(y^{-1}x)^2 =b^{-1}ab^{-1}a
=a^2$, so $(y^{-1}x)^6 = x^2$.  We have $y^{-1}xy = a^2b^{-1}$,and so $y^{-1}xyx^{-1} =b$ and $(y^{-1}xyx^{-1})^3 = 1$. Also $a = a^3 a^{-2} =
x(y^{-1}x)^{-2} = yx^{-1}y$.  Note  - we use it later - that $(xy)^6  = (y^{-1}x^{-1})^{-6} = (y^{-1}xx^{-2})^{-6} = (y^{-1}x)^{-6}x^{12} = x^{10}$.

The group  $G$ is generated by four   elements $a,  b, c$ and $d$,
subject to an infinite set of defining relations as follows.   Firstly  $c^{-1}dc = d^2$, so that $c, d$ generate a subgroup $B$ isomorphic to 
the soluble Baumslag-Solitar group $BS(1,2)$.  Also $a^3 = d$, together with the relations of $A$,
$b^3 =1, a^{-1}ba = b^{-1}$.
The remaining relations are defined inductively.   Put $d = d_1, a = a_1$ and  $d_{i+1} =cd_ic^{-1}$ so that $d_{i+1}^2 = d_i$.
Put $d_0 = d_1^2$ and $a_0 = a^2b$.  Then, as above, the subgroup $A = \langle a, b\rangle =\langle a_0, d_1\rangle$.   Now define inductively
$a_{i+1} = a_id_{i+1}^{-1}a_i, b_{i+1} = a_i^{-1}d_{i+1}a_id_{i+1}^{-1} $ and add the relations $b_{i+1}^3 =1, a_{i+1}^{-1}b_{i+1}a_{i+1} = b_{i+1}^{-1}$ for each $i$
to make $A_{i+1} = \langle a_{i+1}, b_{i+1}\rangle \cong A$.  Note that for $i =1$ we have $a = a_1 = a_0d_1^{-1}a_0 = yx^{-1}y$ as above.
The group $G$ is best understood in terms of the subgroup lattice shown in Fig \ref{Fig 2} and the folding sequence shown in Fig \ref{Fig 3}.
Folding operations are described in  \cite{[DJ1]}.   The sequence here only involves Type II folds and vertex morphisms.
In a Type II fold, edges in the same orbit are folded together.   The stabilizer of a representative edge in the orbit is increased from
$E$ to $\langle E, g\rangle $, and the stabilizer of the orbit of the terminal vertex is increased from $U$ to $\langle U, g\rangle $.  Here
$g$ is an element of the representative vertex group $V$ of the initial vertex.   It is possible that the initial vertex and terminal vertex are 
in the same orbit, i.e. $U = V$.   A vertex morphism involves a homomorphism of a particular vertex group that restricts to an isomorphism
on any incident edge group.  Such a homomorphism induces a morphism of the trees associated with the graph of groups and
a homomorphism of the corresponding fundamental groups.   These morphisms are described in detail in \cite{[D00]}.
In fact we do not use vertex morphisms in our construction as explained below.

\begin{figure}

\centering
\begin{tikzpicture}[scale=2]

    \path (0,0) coordinate (p1);
    \path (0, 2) coordinate (p2);
    \path (3, -1) coordinate (p3);
    \path (2, -2) coordinate (p4);
    \path (1,1) coordinate (p5);
    \path (1,-1) coordinate (p6);
   \filldraw (p1) circle (1pt) ;
   
     \filldraw (p2) circle (1pt);
  \filldraw (p3) circle (1pt);  
 
  \filldraw (p4) circle (1pt);
    \filldraw (p5) circle (1pt) ;

    \path (1,3) coordinate (p7);
    \path (4,0) coordinate (p8);
    \path (2,2) coordinate (p9);
    \path (5,1) coordinate (p10);

     \filldraw (p7) circle (1pt);
  \filldraw (p8) circle (1pt);  
 
  \filldraw (p9) circle (1pt);
    \filldraw (p10) circle (1pt) ;

     \draw (p1) -- (p2) -- (p3) --(p4) -- (p1) ;
     \draw (p3)--(6,2) ;
       \draw (p10)--(4,2) ;
         \draw (p9)--(2,2.5) ;
         \draw (p5)--(p7)--(p8) ;
 
 \draw (2.6, -1.6) node [right] {$2$} ;
\draw (-.1,-.1) node [left] {$a^2b = a_0$} ;
\draw (0, 2.1) node [above] {$A=A_1$} ;
\draw (1,3.15) node [above] {$A_2$} ;

\draw (3.1, -1) node [right] {$a^3 = d_1$} ;
\draw (2,-2) node [below] {$a^6 = d_0$} ;
\draw (1.1,1.1) node [right] {$a = a_1$} ;

\draw (2.1,2.1) node [right] {$a_2$} ;

\draw (3.6, -.6) node [right] {$2$} ;
\draw (4.6,.4) node [right] {$2$} ;

\end{tikzpicture}
\vskip 1cm \caption{Subgroup lattice in $G$}\label{Fig 2}
\end{figure}
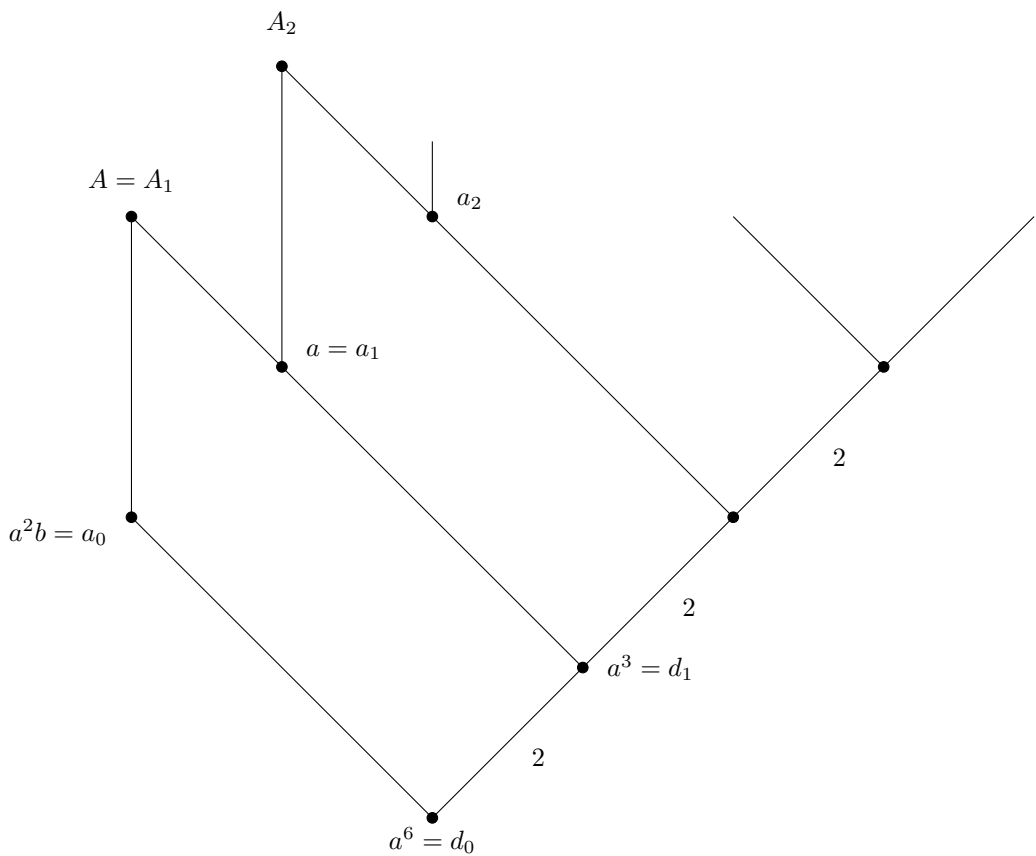
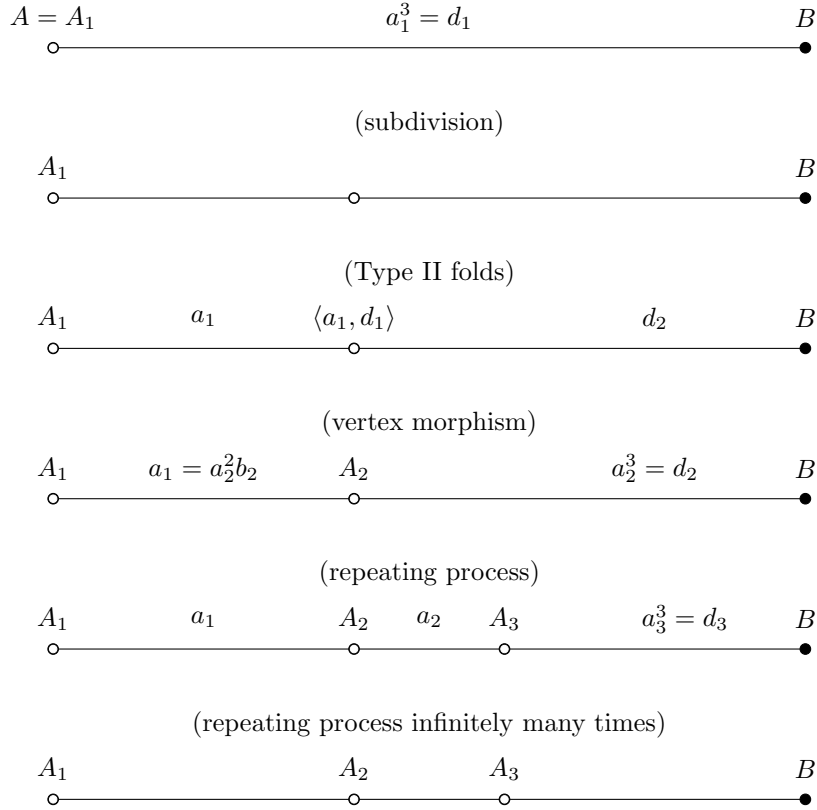
\begin{figure}

\centering
\begin{tikzpicture}[scale=2]

    \path (0,0) coordinate (p1);
    \path (0,1) coordinate (p2);
    \path (0,2) coordinate (p3);
    \path (0,3) coordinate (p4);
    \path (0,4) coordinate (p5);
    \path (0,5) coordinate (p6);
    
    \path (5,0) coordinate (q0);  
    \path (5,1) coordinate (q1);
    \path (5,2) coordinate (q2);
    \path (5,3) coordinate (q3);
    \path (5,4) coordinate (q4);
    \path (5,5) coordinate (q5);
    
    \path (2,0) coordinate (r0);  
     \path (2,1) coordinate (r1);
    \path (2,2) coordinate (r2);
    \path (2,3) coordinate (r3);
    \path (2,4) coordinate (r4);

\path (3,0) coordinate (s0);  
\path (3,1) coordinate (s1);

   \filldraw (p1) circle (1pt) ;
   
     \filldraw (p2) circle (1pt);
  \filldraw (p3) circle (1pt);  
 
  \filldraw (p4) circle (1pt);
    \filldraw (p5) circle (1pt) ;
     \filldraw (p6) circle (1pt) ;
  \filldraw (q0) circle (1pt) ;

     \filldraw (q1) circle (1pt) ;
   
     \filldraw (q2) circle (1pt);
  \filldraw (q3) circle (1pt);  
 
  \filldraw (q4) circle (1pt);
    \filldraw (q5) circle (1pt) ;

      \filldraw (r0) circle (1pt) ;

      \filldraw (r1) circle (1pt) ;
   
     \filldraw (r2) circle (1pt);
  \filldraw (r3) circle (1pt);  
 
  \filldraw (r4) circle (1pt);
  
    \filldraw (s0) circle (1pt) ;

   \filldraw (s1) circle (1pt) ;
\draw (q5)--(p6);
    \draw (q4)--(p5);
    \draw (q3)--(p4);
\draw (q2)--(p3);
\draw (q1)--(p2);
\draw (q0)--(p1);
\draw (2.5, 5.2) node  {$a_1^3 = d_1$};
\draw (2.5,4.5) node  {(subdivision)};
\draw (1,3.2) node  {$a _1$} ;
\draw (4, 3.2) node  {$d_2$} ;
\draw (2.5, 3.5) node  {(Type II folds)};
\draw (2, 3.2) node  {$\langle a_1, d_1\rangle$} ;
\draw (1, 2.2) node  {$a_1 = a_2^2b_2$};
\draw (4, 2.2) node  {$a_2^3 = d_2$} ;
\draw (2.5,2.5) node  {(vertex morphism)};
\draw (2.5, 1.5) node {(repeating process)};
\draw (1, 1.2) node {$a_1$} ;
\draw (2.5, 1.2) node  {$a_2$} ;
\draw (4.2, 1.2) node  {$a_3^3 = d_3$} ;
\draw (2.5, .5) node  {(repeating process infinitely many times)} ;
\draw  (0,.2)  node {$A_1$} ;
\draw  (0,1.2)  node {$A_1$} ;
\draw  (0,2.2)  node {$A_1$} ;
\draw  (0,3.2)  node {$A_1$} ;
\draw  (0,4.2)  node {$A_1$} ;
\draw  (0,5.2)  node {$A=A_1$} ;
\draw  (5,.2)  node {$B$} ;
\draw  (5,1.2)  node {$B$} ;
\draw  (5,2.2)  node {$B$} ;
\draw  (5,3.2)  node {$B$} ;
\draw  (5,5.2)  node {$B$} ;
\draw  (5,4.2)  node {$B$} ;

\draw  (2,.2)  node {$A_2$} ;
\draw  (2,1.2)  node {$A_2$} ;
\draw  (2,2.2)  node {$A_2$} ;

\draw  (3,.2)  node {$A_3$} ;
\draw  (3,1.2)  node {$A_3$} ;

  \filldraw  (p1) [white] circle (.7pt) ;
    \filldraw  (p2) [white] circle (.7pt) ;
 \filldraw  (p3) [white] circle (.7pt) ;
 \filldraw  (p4) [white] circle (.7pt) ;
 \filldraw  (p5) [white] circle (.7pt) ;
 \filldraw  (r1) [white] circle (.7pt) ;
 \filldraw  (r2) [white] circle (.7pt) ;
 \filldraw  (r0) [white] circle (.7pt) ;
 \filldraw  (s0) [white] circle (.7pt) ;
 \filldraw  (s1) [white] circle (.7pt) ;
 
 \filldraw  (p6) [white] circle (.7pt) ;
 \filldraw  (r4) [white] circle (.7pt) ;
 \filldraw  (r3) [white] circle (.7pt) ;

  \end{tikzpicture}

\vskip 1cm \caption{Folding sequence of graph of groups}\label{Fig 3}
\end{figure}

In \cite{[DJ1]} it is shown that $G \cong A *_CG $ where $A = \langle a, b \rangle  = \langle a_0, d_1\rangle$ and $C = \langle a_1\rangle$.
Let $D = \langle d_1 \rangle$.
Let $Y$ be the $G$-graph with two orbits of vertices
$VY = \{ gA, gD| g \in G\} $ and two orbits of edges $EY = \{ (gA, gD), (gD, gcD) | g\in G\}$.

\begin{figure}

\centering
\begin{tikzpicture}[scale=2]

    \path (0,0) coordinate (p1);
    \path (0,1) coordinate (p2);
    \path (0,2) coordinate (p3);
    \path (0,3) coordinate (p4);
    \path (0,4) coordinate (p5);
    
    \path (5,0) coordinate (q0);  
    \path (5,1) coordinate (q1);
    \path (5,2) coordinate (q2);
    \path (5,3) coordinate (q3);
    \path (5,4) coordinate (q4);
    
    \path (2,0) coordinate (r0);  
     \path (2,1) coordinate (r1);
    \path (2,2) coordinate (r2);

\path (3,0) coordinate (s0);  
\path (3,1) coordinate (s1);
   
     \filldraw (p2) circle (1pt);
  \filldraw (p3) circle (1pt);  
 
  \filldraw (p4) circle (1pt);
    \filldraw (p5) circle (1pt) ;
  \filldraw (q0) circle (1pt) ;

     \filldraw (q1) circle (1pt) ;
   
     \filldraw (q2) circle (1pt);
  \filldraw (q3) circle (1pt);  
 
  \filldraw (q4) circle (1pt);

      \filldraw (r0) circle (1pt) ;

      \filldraw (r1) circle (1pt) ;
   
     \filldraw (r2) circle (1pt);
 
  
    \filldraw (s0) circle (1pt) ;

   \filldraw (s1) circle (1pt) ;
    \draw (q4)--(p5);
    \draw (q3)--(p4);
\draw (q2)--(p3);
\draw (q1)--(p2);
\draw (q0)--(p1);
\draw (2.5, 4.2) node  {$a_1^3 = d_1$};
\draw (2.5, 1.5) node  {(subdivision and Type II folds)};
\draw (2.5, 3.5) node  {(Type II fold)};

\draw (2, 2.2) node  {$A_2$} ;
\draw (1, 2.2) node  {$a_1 = a_2^2b_2$};
\draw (1, 1.2) node  {$a_1$};

\draw (4, 2.2) node  {$a_2^3 = d_2$} ;
\draw (2.5,2.5) node  {(subdivision and Type II folds)};
\draw (2.5, 1.2) node  {$a_2$} ;
\draw (4.2, 1.2) node  {$a_3^3 = d_3$} ;
\draw (2.5, .5) node  {(repeating  infinitely many times)} ;
\draw  (0,.2)  node {$A_1$} ;
\draw  (0,1.2)  node {$A_1$} ;
\draw  (0,2.2)  node {$A_1$} ;
\draw  (0,3.2)  node {$A_1$} ;
\draw  (0,4.2)  node {$A=A_1$} ;
\draw  (5,1.2)  node [left] {$d_3$} ;
\draw  (5,2.2)  node [left]{$d_2$} ;
\draw  (5,3.2)  node [left]{$d_2$} ;
\draw  (5,4.2)  node [left] {$d_1$} ;

\draw  (2,.2)  node {$A_2$} ;
\draw  (2,1.2)  node {$A_2$} ;
 \filldraw (p1) circle (1pt);
   \filldraw  (p1) [white] circle (.7pt) ;
    \filldraw  (p2) [white] circle (.7pt) ;
 \filldraw  (p3) [white] circle (.7pt) ;
 \filldraw  (p4) [white] circle (.7pt) ;
 \filldraw  (p5) [white] circle (.7pt) ;
 \filldraw  (r1) [white] circle (.7pt) ;
 \filldraw  (r2) [white] circle (.7pt) ;
 \filldraw  (r0) [white] circle (.7pt) ;
 \filldraw  (s0) [white] circle (.7pt) ;
 \filldraw  (s1) [white] circle (.7pt) ;

\draw  (3,.2)  node {$A_3$} ;
\draw  (3,1.2)  node {$A_3$} ;
\draw (5.4,1 ) [dashed]  circle (.4cm) ;
\draw (5.4,0 ) [dashed]  circle (.4cm) ;
\draw (5.4,2 ) [dashed]  circle (.4cm) ;
\draw (5.4,3 ) [dashed]  circle (.4cm) ;
\draw (5.4,4 ) [dashed]  circle (.4cm) ;
\draw (5.8,0) node [right] {$\langle d_1, d_2, d_3, \dots \rangle $} ;
\draw (5.8,1) node [right] {$d_3$} ;
\draw (5.8,2) node [right] {$d_2$} ;
\draw (5.8,3) node [right] {$d_2$} ;
\draw (5.8,4) node [right] {$d_1$} ;

  \end{tikzpicture}

\vskip 1cm \caption{Folding sequence of graphs}\label{Fig 4}
\end{figure}
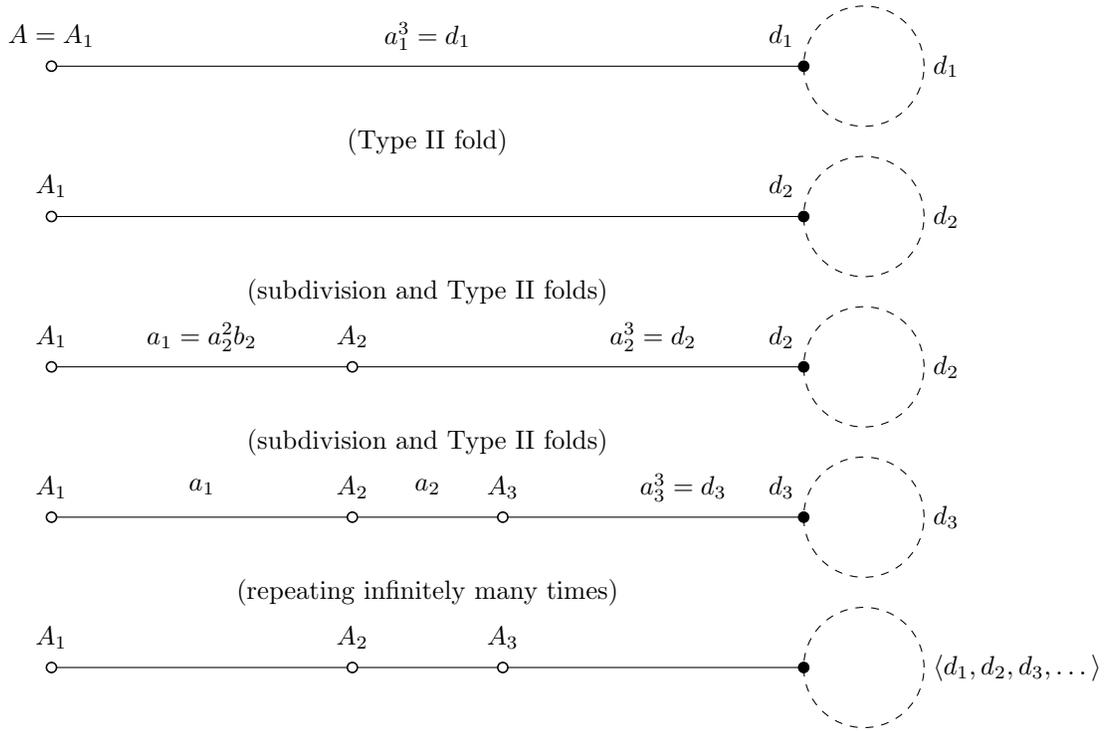

In $Y$ the vertex $A$ is incident with $[A, D] = 9$ edges, as is every vertex in the same orbit.
The vertex $D$ is incident with $4$ edges.   One edge in one edge orbit  connects $D$ to $A$ and  there are
three edges in the other orbit connecting $D$ to $cD, c^{-1}D$ and $dc^{-1}D$.  Note that
$d = d_1$ fixes the edge $(D, cD)$ and transposes the edges $(D, c^{-1}D), (D, dc^{-1}D)$.
If one removes the edges of $Y$ in the first orbit one is left with a set of $3$-regular trees.
If one directs these subgraphs by putting an arrow from $D$ to $cD$, then every vertex has one edge pointing
away from it and two pointing towards it.
The graph $Y$ is connected because $G$ is generated by $A, D$ and $c$.
One obtains a vertex transitive $G$-graph $X$ from $Y$ by taking the  orbit of vertices containing $D$ and joining two
vertices by an edge if they are joined by an edge in $Y$, or they are not joined by an edge in $Y$ but are distance two apart in $Y$.   In $X$ each vertex will have degree
$10$.  Thus $D$ is a vertex in $X$.  It has $2$ vertices adjacent to it which were already adjacent to it in $Y$.  The one vertex in $Y$ adjacent 
to $D$ 
in $Y$ which is not in $X$ has $9$ adjacent vertices including $D$ itself, the $8$ other vertices will be adjacent to $D$ in $X$.   It is easier to work with the graph $Y$, which is quasi-isometric to $X$.
In Fig \ref{Fig 4} a sequence of folding operations is described for the graph $Y$.
These are similar to those of Fig \ref {Fig 3}.   However in Fig \ref {Fig 3} the operations are for trees. Vertex morphisms
are included which change the group acting.    In Fig \ref {Fig 4} the operations are on $G$-graphs in which
the group acting remains the same throughout.  Thus we are assuming that all the vertex morphisms
have been carried out before we start. 
The first diagram in Fig \ref {Fig 4} shows the graph $G\backslash Y$.  Each edge of the quotient graph
is labelled by its stabilizer in a lift to $Y$, as in Bass-Serre theory.
The first folding operations results in the edges at $D$  (labelled with a $\bullet$ in Fig \ref {Fig 4}) in the same $d$ orbit being folded together.
The stabilizer of $D$ is increased, as are the stabilizers of all the edges in the orbit of $(D, cD)$. A new stabilizer  includes
the original stabilizer as a subgroup of index two.
The degree of $D$ changes to $5$ as  $D$ is identified with $d_2D$ and the two edges $(A, D), (d_2A, D)$ are now incident
with the new vertex. The graph still contains $3$-regular trees as before.

The next operation, which is subdivision, inserts a $\circ $ vertex  on each edge of the orbit containing $(A, D)$.
The next operation comprises folding from the $\circ$ vertex  labelled $A_1$ and the $\bullet $ vertex labelled $D$ towards the new $\circ $  vertex
just created.   In the folding from $A$, edges corresponding to the three cosets of $D$ which belong to $\langle a \rangle$ 
are folded together.   The vertex $A$ will then have degree $3$.    In the folding from $D$, the edges in the same $\langle d_2\rangle $
are folded together, so that the degree of $D$ again becomes $4$.
The defining relations for $G$ ensure that the created vertex has  stabilizer $A_2$.
In the graph now obtained, the vertices in the orbit of  $A = A_1$ are cut points.  Removing one of these vertices
gives three components.    Removing all the vertices in this orbit with all their incident edges gives an infinite set of component 
graphs, each of which is isomorphic to $Y$.  Thus one can repeat this process on each of these graphs as indicated in Fig 4.

In fact in $Y$, before any folding operation, each of the $\circ $ vertices is a cut point.   There is a tree decomposition of $Y$ as in \cite{[DK]},
in which the $G$-tree $T$ has two orbits of vertices.    One orbit corresponds to the  $2$-blocks, where each $2$-block corresponds to 
a maximal $2$-connected subgraph, and the other orbit corresponds to the cut points.

Removing a particular $\circ$ vertex from $Y$ results in $3$ components.  In each  $2$-block a $\circ$-vertex has degree $3$.

If we consider a cycle in $Y$, then under the successive subdivision and folding operations,  the cycle will eventually be 
a subtree.  But there will have to be at least one fold at a vertex  at each stage in the iteration.   Thus if we start with a cycle with $k$ edges
then after $k$ iterations, the image of the cycle will be a subtree. We give a more precise explanation of how this happens after considering an example.

\begin{figure}

\centering
\begin{tikzpicture}[scale=.8]
   \path (3,0) coordinate (p1);
    \path (5,0) coordinate (p2);
    \path (.5, 3) coordinate (p3);
    \path (1.5,1) coordinate (p4);
    \path (7.5,3) coordinate (p5);
    \path (6.5,1) coordinate (p6);

    \path (3,8) coordinate (q1);
    \path (5,8) coordinate (q2);
    \path (.5, 5) coordinate (q3);
    \path (1.5,7) coordinate (q4);
    \path (7.5,5) coordinate (q5);
    \path (6.5,7) coordinate (q6);

    \path (-1,4) coordinate (r1);
    \path (9,4) coordinate (r2);
    \path (1.5,9) coordinate (r3);
    \path (6.5,9) coordinate (r4);
    \path (1.5,-1) coordinate (r5);
    \path (6.5,-1) coordinate (r6);
    
  \path (4,0) coordinate (s1);
    \path (1,2) coordinate (s2);
    \path (1,6) coordinate (s3);
    \path (7,6) coordinate (s4);
    \path (7,2) coordinate (s5);
    \path (4,8) coordinate (s6);

       \filldraw (p1) circle (3pt) ;
   
     \filldraw (p2) circle (3pt);
  \filldraw (p3) circle (3pt);  
 
  \filldraw (p4) circle (3pt);
    \filldraw (p5) circle (3pt) ;
     \filldraw (p6) circle (3pt) ;

       \filldraw (q1) circle (3pt) ;
   
     \filldraw (q2) circle (3pt);
  \filldraw (q3) circle (3pt);  
 
  \filldraw (q4) circle (3pt);
    \filldraw (q5) circle (3pt) ;
     \filldraw (q6) circle (3pt) ;

       \filldraw (r1) circle (3pt) ;
   
     \filldraw (r2) circle (3pt);
  \filldraw (r3) circle (3pt);  
 
  \filldraw (r4) circle (3pt);
    \filldraw (r5) circle (3pt) ;
     \filldraw (r6) circle (3pt) ;

  \filldraw (s1) circle (1.5pt) ;
   
     \filldraw (s2) circle (1.5pt);
  \filldraw (s3) circle (1.5pt);  
 
  \filldraw (s4) circle (1.5pt);
    \filldraw (s5) circle (1.5pt) ;
     \filldraw (s6) circle (1.5pt) ;

     \draw (p1) -- (p2) ;
      \draw (p3) -- (p4);
 \draw (p5) -- (p6);

     \draw (q1) -- (q2) ;
      \draw (q3) -- (q4);
 \draw (q5) -- (q6);

    \draw [dashed] (p3)--(r1)--(q3) ;
    \draw [dashed] (p4)--(r5)--(p1);
      \draw [dashed] (p6)--(r6)--(p2) ;
    \draw [dashed] (q6)--(r4)--(q2);
  \draw [dashed] (q4)--(r3)--(q1) ;
    \draw [dashed] (p5)--(r2)--(q5);

  \path (3,-8) coordinate (p1);
    \path (4, -7.6) coordinate (p2);
    \path (5, -8) coordinate (p3);
    \path (5.4, -6.9) coordinate (p4);
    \path (6,-6) coordinate (p5);
    \path (5.4,-5.1) coordinate (p6);

    \path (5,-4) coordinate (q1);
    \path (4,-4.4) coordinate (q2);
    \path (3,-4) coordinate (q3);
    \path (2.6,-5.1) coordinate (q4);
    \path (2,-6) coordinate (q5);
    \path (2.6,-6.9) coordinate (q6);
    
    \draw (p1)--(p2)--(p3)--(p4)--(p5)--(p6)--(q1)--(q2)--(q3)--(q4)--(q5)--(q6)--(p1) ;
    
   \path (2.5, -9) coordinate (r1);
    \path (5.5, -9) coordinate (r2);
    \path (1,-6) coordinate (r3);
    \path (7,-6) coordinate (r4);
    \path (2.5,-3) coordinate (r5);
    \path (5.5,-3) coordinate (r6);
    
    \draw [dashed] (r1)--(p1) ;
    \draw [dashed] (r2)--(p3);
    \draw [dashed] (r3)--(q5);
    \draw [dashed] (r4)--(p5);
    \draw [dashed] (r5)--(q3);
    \draw [dashed] (r6)--(q1);
     
       \filldraw (p1) circle (3pt) ;

    \filldraw (p3) circle (3pt);

    \filldraw (p5) circle (3pt) ;

       \filldraw (q1) circle (3pt) ;

     \filldraw (q3) circle (3pt);

    \filldraw (q5) circle (3pt) ;
    
       \filldraw (q2)  circle (1.5pt) ;

   \filldraw (q4) circle (1.5pt);
  \filldraw (q6) circle (1.5pt);  
 
  \filldraw  (p2) circle (1.5pt);
    \filldraw  (p4) circle (1.5pt) ;
     \filldraw (p6) circle (1.5pt) ;

\draw (4, .3) node  {$a$} ;
\draw (1.3, 6) node  {$a$} ;
\draw (6.7, 5.9) node  {$a$} ;
\draw (4, 7.7) node  {$a$} ;
\draw (1.3,2.1) node  {$a$} ;
\draw (6.7,2.1) node  {$a$} ;
\draw (5.6,-6) node  {$d_2$} ;
\draw (2.4,-6) node  {$d_2$} ;
\draw (3.1,-7.6) node  {$d_2$} ;
\draw (4.9, -7.6) node  {$d_2$} ;
\draw (3.1, -4.4) node  {$d_2$} ;
\draw (4.9,-4.4) node  {$d_2$} ;

       \filldraw (r1) circle (3pt) ;
   
     \filldraw (r2) circle (3pt);
  \filldraw (r3) circle (3pt);  
 
  \filldraw (r4) circle (3pt);
    \filldraw (r5) circle (3pt) ;
     \filldraw (r6) circle (3pt) ;

  \end{tikzpicture}

\vskip 1cm \caption{Folding a cycle I}\label{Fig 5}
\end{figure}
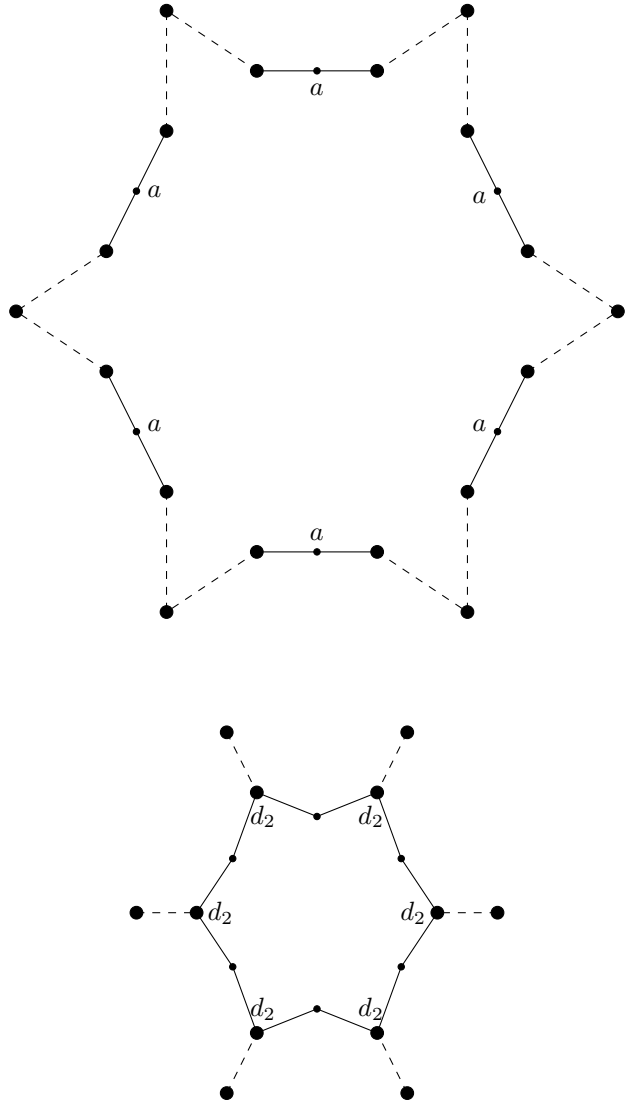

\begin{figure}

\centering
\begin{tikzpicture}[scale=1.5]

  \path (3,1) coordinate (p1);
    \path (4, 1.4) coordinate (p2);
    \path (5, 1) coordinate (p3);
    \path (5.4, 2.1) coordinate (p4);
    \path (6,3) coordinate (p5);
    \path (5.4,3.9) coordinate (p6);

    \path (5,5) coordinate (q1);
    \path (4,4.6) coordinate (q2);
    \path (3,5) coordinate (q3);
    \path (2.6,3.9) coordinate (q4);
    \path (2,3) coordinate (q5);
    \path (2.6,2.1) coordinate (q6);
    
    \draw (p1)--(p2)--(p3)--(p4)--(p5)--(p6)--(q1)--(q2)--(q3)--(q4)--(q5)--(q6)--(p1) ;
    
   \path (4, -2) coordinate (r1);
    \path (4,-4) coordinate (r2);
    \path (3.25, -2.5) coordinate (r3);
    \path (4.75, -3.5) coordinate (r4);
    \path (3.25, -3.5) coordinate (r5);
    \path (4.75, -2.5) coordinate (r6);

  \path (3.5, -2) coordinate (s1);
    \path (4.5, -4) coordinate (s2);
    \path (3,-3) coordinate (s3);
    \path (5,-3) coordinate (s4);
    \path (4.5, -2) coordinate (s5);
    \path (3.5, -4) coordinate (s6);

    \draw (r1)--(r2);
     \draw (r3)--(r4);
      \draw (r5)--(r6);
       \draw (s1)--(s2);
        \draw (s3)--(s4);
         \draw (s5)--(s6);

  \path (3, -1) coordinate (t1);
    \path (5,-1) coordinate (t2);
    \path (2,-3) coordinate (t3);
    \path (3,-5) coordinate (t4);
    \path (5,-5) coordinate (t5);
    \path (6,-3) coordinate (t6);

    \draw [dashed] (t1)--(s1);
     \draw  [dashed] (t2)--(s5);
      \draw  [dashed] (t3)--(s3);
       \draw [dashed]  (t4)--(s6);
        \draw [dashed]  (t5)--(s2);
         \draw [dashed]  (t6)--(s4);

  \path (2.5,0) coordinate (w1);
    \path (5.5,0) coordinate (w2);
    \path (1,3) coordinate (w3);
    \path (7,3) coordinate (w4);
    \path (2.5,6) coordinate (w5);
    \path (5.5,6) coordinate (w6);

    \draw [dashed] (w1)--(p1);
     \draw  [dashed] (w2)--(p3);
      \draw  [dashed] (w3)--(q5);
       \draw [dashed]  (w4)--(p5);
        \draw [dashed]  (w5)--(q3);
         \draw [dashed]  (w6)--(q1);
         
              \filldraw (p1) circle (2pt) ;
   
     \filldraw (p2) circle (1pt);
  \filldraw (p3) circle (2pt);  
 
  \filldraw (p4) circle (1pt);
    \filldraw (p5) circle (2pt) ;
     \filldraw (p6) circle (1pt) ;

       \filldraw (q1) circle (2pt) ;
   
     \filldraw (q2) circle (1pt);
  \filldraw (q3) circle (2pt);  
 
  \filldraw (q4) circle (1pt);
    \filldraw (q5) circle (2pt) ;
     \filldraw (q6) circle (1pt) ;
        
              \filldraw (t1) circle (2pt) ;
   
     \filldraw (t2) circle (2pt);
  \filldraw (t3) circle (2pt);  
 
  \filldraw (t4) circle (2pt);
    \filldraw (t5) circle (2pt) ;
     \filldraw (t6) circle (2pt) ;

              \filldraw (w1) circle (2pt) ;
   
     \filldraw (w2) circle (2pt);
  \filldraw (w3) circle (2pt);  
 
  \filldraw (w4) circle (2pt);
    \filldraw (w5) circle (2pt) ;
     \filldraw (w6) circle (2pt) ;

       \filldraw (r1) circle (1pt) ;
   
     \filldraw (r2) circle (1pt);
  \filldraw (r3) circle (1pt);  
 
  \filldraw (r4) circle (1pt);
    \filldraw (r5) circle (1pt) ;
     \filldraw (r6) circle (1pt) ;

       \filldraw (s1) circle (2pt) ;
   
     \filldraw (s2) circle (2pt);
  \filldraw (s3) circle (2pt);  
 
  \filldraw (s4) circle (2pt);
    \filldraw (s5) circle (2pt) ;
     \filldraw (s6) circle (2pt) ;

 \filldraw (2.3,  3.45) circle (1.5pt);
  \filldraw (2.3,  3.45) [white] circle (1pt);
   \filldraw (2.8,  1.55) circle (1.5pt);
  \filldraw (2.8,  1.55) [white] circle (1pt);
 \filldraw (5.7,  3.45) circle (1.5pt);
  \filldraw (5.7,  3.45) [white] circle (1pt);
 \filldraw (5.7, 2.55) circle (1.5pt);
  \filldraw (5.7,2.55) [white] circle (1pt);
 \filldraw (3.5,1.2) circle (1.5pt);
  \filldraw (3.5,1.2) [white] circle (1pt);
 \filldraw (4.5,1.2) circle (1.5pt);
  \filldraw (4.5,1.2) [white] circle (1pt);
   \filldraw (5.2,  1.55) circle (1.5pt);
  \filldraw (5.2, 1.55) [white] circle (1pt);
 \filldraw (5.7,3.45) circle (1.5pt);
  \filldraw (5.7,  3.45) [white] circle (1pt);
 \filldraw (5.2,4.45) circle (1.5pt);
  \filldraw (5.2,4.45) [white] circle (1pt);
 \filldraw (2.3, 2.55) circle (1.5pt);
  \filldraw (2.3, 2.55) [white] circle (1pt);
 \filldraw (4,-3) circle (1.5pt);
  \filldraw (4,-3) [white] circle (1pt);

 \filldraw (4.5,4.8) circle (1.5pt);
  \filldraw (4.5,4.8) [white] circle (1pt);
 \filldraw (3.5,4.8) circle (1.5pt);
  \filldraw (3.5,4.8) [white] circle (1pt);
 \filldraw (2.8,4.45) circle (1.5pt);
  \filldraw (2.8,4.45) [white] circle (1pt);

\draw (4,4.4) node  {$a$} ;
\draw (5.3,3.7) node  {$a$} ;
\draw (2.7,2.2) node  {$a$} ;
\draw  (5.3,2.2)  node {$a$} ;
\draw  (4,1.6)  node {$a$} ;
\draw  (2.7,3.7)  node {$a$} ;

\draw (5.7,3) node  {$d_2$} ;
\draw (2.3,3) node  {$d_2$} ;
\draw (3.1,1.3) node  {$d_2$} ;
\draw  (4.9,1.3)  node {$d_2$} ;
\draw  (3.1,4.7)  node {$d_2$} ;
\draw  (4.9,4.7)  node {$d_2$} ;

\draw (3.4,-2) [left] node  {$x$} ;
\draw (4.6,-2) [right] node  {$x$} ;
\draw (3.4,-4)[left] node  {$x$} ;
\draw  (4.6,-4) [right] node {$x$} ;
\draw  (5,-2.85)  [right] node {$x$} ;
\draw  (3,-2.85) [left]  node {$x$} ;

\draw (3.15,-2.5) [left] node  {$y$} ;
\draw (4.85,-2.5) [right] node  {$y$} ;
\draw (4,-2)[above] node  {$y$} ;
\draw  (4,-4) [below] node {$y$} ;
\draw  (3.15, -3.5)  [left] node {$y$} ;
\draw  (4.85, -3.5) [right]  node {$y$} ;

  \end{tikzpicture}

\vskip 1cm \caption{Folding a cycle II}\label{Fig 6}
\end{figure}

To illustrate the remarks above, in Fig \ref {Fig 5}  and Fig \ref {Fig 6}  the effect of the folding sequence
is shown on a particular cycle in $Y$.  This particular cycle is reduced to a subtree after one iteration of the folding sequence.
Probably it is a shortest cycle in $Y$.
Vertices in the orbit of $A$ in $Y$ are indicated by a small $\bullet  $, the other vertices are indicated by a larger $\bullet $.
Vertices created in the process by subdivision are indicated with a $\circ $.
There are two orbits of edges.   The ones in the orbit of edges incident with an $A$-orbit  vertex are indicated with
a continuous line.  These are called {\it solid} edges.  The others are indicated with a dashed line are called {\it dashed} edges.
As the cycle passes through an $A$-orbit  vertex,  the different directions one can proceed correspond to the
nine cosets of $\langle a^3\rangle $ in $A$.  As folding takes place at each vertex in the cycle, the direction taken must correspond to
one of the two cosets containing either $a$ or $a^2$.  The diagram shows the choice at each such vertex. (We always make the same choice $a$.)
At a $\bullet$ vertex, if one is proceeding from a solid edge to a dashed edge, one proceeds along the only edge directed away
from the vertex. (Recall that every $\bullet $ vertex has one dashed edge directed away from it and two directed towards it.)
If one arrives at a $\bullet$ vertex along a dashed edge and leaves along a dashed edge, then one leaves along the other dashed edge
directed towards the vertex.  The first diagram indicates a uniqe path in $Y$ which in fact turns out to be a cycle.
The fact that one has a cycle is because $(xy)^6 = x^{10} \in D$  fixes an edge of $Y$.

\begin{figure}

\centering
\begin{tikzpicture}[scale=1.1]
   \path (3,0) coordinate (p1);
    \path (5,0) coordinate (p2);
    \path (.5, 3) coordinate (p3);
    \path (1.5,1) coordinate (p4);
    \path (7.5,3) coordinate (p5);
    \path (6.5,1) coordinate (p6);
    
      \path (0, 3.5) coordinate (a1);
    \path (8,3.5) coordinate (a2);
    \path (1.6, 7.7) coordinate (a3);
    \path (6.4,7.7) coordinate (a4);
    \path (1.6, .3) coordinate (a5);
    \path (6.4, .3) coordinate (a6);
    
      \path (0, 4.5) coordinate (b1);
    \path (8, 4.5) coordinate (b2);
    \path (2.4, 8.2) coordinate (b3);
    \path (5.6, 8.2) coordinate (b4);
    \path (2.4, -.2) coordinate (b5);
    \path (5.6, -.2) coordinate (b6);

  \path (4, -3) coordinate (c1);
    \path (4,11) coordinate (c2);
    \path (-2, 0) coordinate (c3);
    \path (-2, 8) coordinate (c4);
    \path (10, 8) coordinate (c5);
    \path (10, 0) coordinate (c6);

  \path (3.5, -1) coordinate (d1);
    \path (4.5, -1) coordinate (d2);
    \path (0, 2) coordinate (d3);
    \path (.45, 1.25) coordinate (d4);
    \path (8, 2) coordinate (d5);
    \path (7.55, 1.25) coordinate (d6);

      \path (3.75, -2) coordinate (e1);
    \path (4.25, -2) coordinate (e2);
    \path (-1,1) coordinate (e3);
    \path (-.78, .625) coordinate (e4);
    \path (9, 1) coordinate (e5);
    \path (8.78, .625) coordinate (e6);

    \path (3,8) coordinate (q1);
    \path (5,8) coordinate (q2);
    \path (.5, 5) coordinate (q3);
    \path (1.5,7) coordinate (q4);
    \path (7.5,5) coordinate (q5);
    \path (6.5,7) coordinate (q6);

    \path (-1,4) coordinate (r1);
    \path (9,4) coordinate (r2);
    \path (1.5,9) coordinate (r3);
    \path (6.5,9) coordinate (r4);
    \path (1.5,-1) coordinate (r5);
    \path (6.5,-1) coordinate (r6);
    
  \path (4,0) coordinate (s1);
    \path (1,2) coordinate (s2);
    \path (1,6) coordinate (s3);
    \path (7,6) coordinate (s4);
    \path (7,2) coordinate (s5);
    \path (4,8) coordinate (s6);

       \filldraw (p1) circle (2pt) ;
   
       \filldraw (a1) circle (2pt) ;
         \filldraw (a2) circle (2pt);
  \filldraw (a3) circle (2pt);  
 
  \filldraw (a4) circle (2pt);
    \filldraw (a5) circle (2pt) ;
     \filldraw (a6) circle (2pt) ;
     
       \filldraw (b1) circle (2pt) ;

  \filldraw (b2) circle (2pt);
  \filldraw (b3) circle (2pt);  
 
  \filldraw (b4) circle (2pt);
    \filldraw (b5) circle (2pt) ;
     \filldraw (b6) circle (2pt) ;

     \filldraw (p2) circle (2pt);
  \filldraw (p3) circle (2pt);  
 
  \filldraw (p4) circle (2pt);
    \filldraw (p5) circle (2pt) ;
     \filldraw (p6) circle (2pt) ;

       \filldraw (q1) circle (2pt) ;
   
     \filldraw (q2) circle (2pt);
  \filldraw (q3) circle (2pt);  
 
  \filldraw (q4) circle (2pt);
    \filldraw (q5) circle (2pt) ;
     \filldraw (q6) circle (2pt) ;

    \filldraw (c1) circle (2pt) ;
   
     \filldraw (c2) circle (2pt);
  \filldraw (c3) circle (2pt);  
 
  \filldraw (c4) circle (2pt);
    \filldraw (c5) circle (2pt) ;
     \filldraw (c6) circle (2pt) ;
   
    \filldraw (d1) circle (1pt) ;
   
     \filldraw (d2) circle (1pt);
  \filldraw (d3) circle (1pt);  
 
  \filldraw (d4) circle (1pt);
    \filldraw (d5) circle (1pt) ;
     \filldraw (d6) circle (1pt) ;

\filldraw (e1) circle (2pt) ;
   
     \filldraw (e2) circle (2pt);
  \filldraw (e3) circle (2pt);  
 
  \filldraw (e4) circle (2pt);
    \filldraw (e5) circle (2pt) ;
     \filldraw (e6) circle (2pt) ;

       \filldraw (r1) circle (2pt) ;
   
     \filldraw (r2) circle (2pt);
  \filldraw (r3) circle (2pt);  
 
  \filldraw (r4) circle (2pt);
    \filldraw (r5) circle (2pt) ;
     \filldraw (r6) circle (2pt) ;

     \draw (p1) -- (d1) --(e1);
      \draw (p2) -- (d2) --(e2) ;
     \draw (p3) -- (d3) --(e3) ;
    \draw (p4) -- (d4) --(e4) ;
  \draw (p5) -- (d5) --(e5) ;
  \draw (p6) -- (d6) --(e6) ;

    \draw [dashed] (p3)--(a1)--(r1)--(b1)--(q3) ;
    \draw [dashed] (p4)--(a5)--(r5)--(b5)--(p1);
      \draw [dashed] (p6)--(a6)--(r6)--(b6)--(p2) ;
    \draw [dashed] (q6)--(a4)--(r4)--(b4)--(q2);
  \draw [dashed] (q4)--(a3)--(r3)--(b3)--(q1) ;
    \draw [dashed] (p5)--(a2)--(r2)--(b2)--(q5);
\draw [dashed] (e1) --(c1)--(e2) ;
 \draw [dashed] (e3) --(c3)--(e4) ;
  \draw [dashed] (e5) --(c6)--(e6) ;

  \path (3.5, 9) coordinate (d1);
    \path (4.5, 9) coordinate (d2);
    \path (0, 6) coordinate (d3);
    \path (.45, 6.75) coordinate (d4);
    \path (8, 6) coordinate (d5);
    \path (7.55, 6.75) coordinate (d6);

      \path (3.75, 10) coordinate (e1);
    \path (4.25, 10) coordinate (e2);
    \path (-1,7) coordinate (e3);
    \path (-.78, 7.375) coordinate (e4);
    \path (9, 7) coordinate (e5);
    \path (8.78, 7.375) coordinate (e6);
\draw [dashed] (e1) --(c2)--(e2) ;
 \draw [dashed] (e3) --(c4)--(e4) ;
  \draw [dashed] (e5) --(c5)--(e6) ;

 \draw (q1) -- (d1) --(e1);
      \draw (q2) -- (d2) --(e2) ;
     \draw (q3) -- (d3) --(e3) ;
    \draw (q4) -- (d4) --(e4) ;
  \draw (q5) -- (d5) --(e5) ;
  \draw (q6) -- (d6) --(e6) ;

    \filldraw (d1) circle (1pt) ;
   
     \filldraw (d2) circle (1pt);
  \filldraw (d3) circle (1pt);  
 
  \filldraw (d4) circle (1pt);
    \filldraw (d5) circle (1pt) ;
     \filldraw (d6) circle (1pt) ;

\filldraw (e1) circle (2pt) ;
   
     \filldraw (e2) circle (2pt);
  \filldraw (e3) circle (2pt);  
 
  \filldraw (e4) circle (2pt);
    \filldraw (e5) circle (2pt) ;
     \filldraw (e6) circle (2pt) ;

   \end{tikzpicture}

\vskip 1cm \caption{A $60$-cycle }\label{Fig 7}
\end{figure}

In general if one starts with a cycle in $Y$, then  after one stage of the iteration the cycle will have become a closed path
and folding will have taken place at at least one vertex.  If the image is not already a subtree (as in the example above) then
further folding must take place at the next stage.   This folding must take place at a point that is at the end of a fold of the previous stage.
Thus it is either at a $\bullet $ vertex which is at the end of two dashed edges which have been folded together and the new fold will also be
between two dashed edges, or it will be at a $\circ $ and it will be between edges which come from two distinct folds at $\bullet $ vertices.
It can be seen that the number of points where folding can take place is strictly less than at the previous stage.
Thus if there are original cycle has $k$ edges (or vertices) then its image is a subtree after $k$ stages.  In fact it will become a subtree
 after many less stages.  In Fig \ref {Fig 7} we show a $60$-cycle in $Y$ that reduces to a tree after two stages of the iteration.  After one stage it will
 be the like the first cycle of Fig 5 with spikes attached.  In Fig \ref{Fig 7}  we preserve the previous convention that edges in the $3$-regular subtrees
 are dashed, and the other edges are shown with continuous lines.   There is a sequence $C_n$ of cycles of increasing size such that
 $C_n$ folds to $C_{n-1}$ with spikes after one iteration of the folding sequence.    If $c_n$ is the number of vertices of $C_n$, then $c_1 =24,
 c_2 = 60, c_3 =132$.   Each of these cycles is such that the distance between two vertices in the cycle is
 the same as that in $Y$.

Let $Z$ be a particular subgraph which is $3$-regular tree consisting of dashed edges.
The way the edges are oriented gives a height function $\phi : VZ \rightarrow \Z$ by defining $\phi (v_0) = 0$ for some fixed vertex $v_0 \in VZ$,
and such that if $e$ is an oriented edge of $Z$ with initial vertex $\iota e$ and terminal vertex $\tau e$, then $\partial \phi (e) = \phi (\tau e) - \phi (\iota  e) = 1$.  Two vertices of $Z$ are joined by a path in $Y$ in which the only vertices in $Z$ are the end vertices if and only if the two vertices have
the same height.  The shortest such path will be much longer than the shortest path joining them in $Z$.  Thus, from Fig \ref {Fig 5},  two vertices at the same height
in $Z$ that are distance two apart, are joined by a path in $Y$ internally disjoint from $Z$ of length $22$ .  And, from Fig \ref {Fig 7} two vertices at the same
height in $Z$ that are distance $4$ apart, are joined by a path in $Y$ of length $56$ which is internally disjoint from $Z$.   The fact that any
two vertices at the same height are joined by a path outside $Z$ means that any two rays in $Z$ represent the same end.

Consider the effect of a quasi-isometry on a graph $U$.   Let $\theta : U \rightarrow W$ be a quasi-isometry.
Then $\theta $ induces a bijection $\ce ( \theta ): \ce(U) \rightarrow \ce (V)$ and $\ce $ takes special ends to special ends.
This is because if $\omega _1$ and $\omega _2$ are ends of $U$ that are  separated by a set of $s$ vertices then $\ce (\theta )(\omega _1)$
and $\ce (\theta ) (\omega _2)$ are separated by a set of $f(s)$ vertices where $f$ is a function of the form $f(x) = cx + d$.  Thus 
 $k(\omega ) =\infty $ if and only if $k(\ce (\theta) (\omega)) = \infty$.  
 
\begin{figure}

\centering
\begin{tikzpicture}[scale=2]

    \path (0,0) coordinate (p1);

    \path (5,0) coordinate (q0);  
   
        \path (4.5,0) coordinate (q1);  
  \path (3.5,0) coordinate (q2);  
  \path (3.8,0) coordinate (q3);  

    \path (2,0) coordinate (r0);

\path (3,0) coordinate (s0);

   \filldraw (p1) circle (1pt) ;
   
  \filldraw (q1) circle (1pt) ;
 \filldraw (q2) circle (1pt) ;
 \filldraw (q3) circle (1pt) ;

   \filldraw (q0) circle (1pt) ;

      \filldraw (r0) circle (1pt) ;

    \filldraw (s0) circle (1pt) ;

 \draw (q0)--(p1);

\draw  (0,.2)  node {$\bar A_1$} ;
\draw  (4.5,.2)  node {$\bar A_k$} ;
\draw  (5.3,.2)  node {$\langle \bar A_{k+1}, H\rangle$} ;

\draw  (2,.2)  node {$\bar A_2$} ;

\draw  (3,.2)  node {$\bar A_3$} ;

  \end{tikzpicture}

\vskip 1cm \caption{Structure tree $S_n$}\label{Fig 9}
\end{figure}
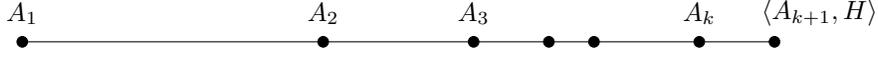

To clarify why the graph $Y$ is not quasi-isometric to a Cayley graph, we construct an inaccessible Cayley graph with similar properties
to $Y$, but point out the significant difference.
We construct an inaccessible group using the lattice of Fig \ref {Fig 2}.    Let $P$ be the subgroup of $G$ generated by all the $A_i$'s,  $i = 1, 2, \dots $.
It can be seen that $P$ is the fundamental group of a graph of groups  $(\cg , N)$  in which the underlying graph $N$ has vertex set which
is the natural numbers $\{1, 2, \dots \}$ and there are edges  $(i, i+1)$ for each $i\in VN$.  The vertex group $\cg (i) = A_i$ and the edge group
corresponding to $(i, i+1)$ is generated by $a_i$.
In $P$ the element $d_0 = a^6$ is central.   If we form the quotient group $\bar P = P/\langle d_0 \rangle $ then $\bar P$ also has a graph of groups 
decomposition with the same underlying graph
and in which the vertex group corresponding to $i$ is $\bar A_1 = A_i/ \langle d_0 \rangle $.
Consider the subgroup $D$  of $P$ generated by $d_1, d_2, \dots $.   This will be locally cyclic.  After factoring out $d_0$ we get
a group $\bar D$ which is locally finite  cyclic.  In fact it is isomorphic to the additive group of dyadic rationals, i.e rationals of the form
$m/2^n$, where $m, n$ are integers.    Note that $P$ is generated by $D$ and $A_1$, so $\bar P$ is generated by $\bar D$ and $\bar A_1$.   Let $H$ be a finitely generated
one-ended  group that contains a subgroup isomorphic to $\bar D$.   Such a group certainly exists.  Any countable group is contained in
a finitely generated group and the direct product of a finitely generated group with a free abelian group of rank two creates a one ended
group.   Form the group $\bar G = \bar P *_{\bar D} H$.   This will be an inaccessible group.
The sequence $S_n$ of structure trees for a Cayley graph   of $\bar G$  will be very similar to the sequence $T_n$ of  structure 
trees for $Y$.  The structure tree $S_n$ is the fundamental group of the graph of groups shown in Fig \ref {Fig 9}. An  edge group of $S_n$ will be a conjugate of the finite cyclic group $\langle a_i\rangle/ \langle d_0 \rangle  $ for some $i$.
Note however that we can choose a generating set for $\bar G$ so that it includes a generating set for $H$, and then the corresponding  Cayley
graph $W$ for $\bar G$ will have a locally finite one-ended subgraph.     This is the important difference with
the graph $Y$. 

The graph $Y$ has countably many subgraphs  which are $3$-regular trees.  These subgraphs are a single orbit under the action of $G$.
Let $Z$ be one of these subgraphs.
Let $W$ be a Cayley graph, and suppose there  are quasi-isometries $\theta : Y \rightarrow W$ and $\phi : W \rightarrow Y$.
Any two rays in $Z$ represent the same end $\omega $, and this end will be special.
Since $G$ is countable, the orbit containing this end is countable.    It is also dense in the subspace of special ends.

Let $\omega '  = \ce (\theta )(\omega )$.  Then $\omega '$ will be a special end of $W$, which is the Cayley graph of an inaccessible group
$Q$.  Since $Q$ is inaccessible, there will be an infinite sequence $Q =Q_1, Q_2, \dots $ where $Q_i$ has a decomposition as a free product
with amalgamation over a finite subgroup in which $Q_{i+1}$ is one of the factors, or $Q_i$ is an HNN-group with vertex group $Q_{i+1}$
and finite edge group.  At least one factor in each decomposition is inaccessible.   If there is an infinite sequence of factorizations in which
there is more than one inaccessible factor infinitely many times (as can in fact happen in some inaccessible groups) then there will be no countable orbit of special ends 
that is dense in the space of all special ends.   Thus for any sequence of decompositions of factors of $Q$ we will eventually obtain 
a term $Q_j$  that for each $i>j $ we have that $Q_i = Q_{i+1}*_{F_{i+1}}Q_{i+1}'$  and $Q_{i+1}'$ is accessible.  In fact if $Q_{i+1}'$ has
an infinite one-ended factor, then $Q$ would contain a thick end $\omega _1$ with  $k(\omega _1)$ finite.   But $Y$ contains no such thick end
and so $Q$ has no such end. 
We are then, very much, as  in the situation of the example above, in which all the $Q_{i+1}'$ factors are finite.
Let $Q'$ be the subgroup of $Q$ generated by all the $Q_i'$'s.  This will have a graph of groups decomposition with infinitely many factors,
in which the $Q_i'$'s are the vertex groups.    This group is not  finitely generated.

Now put $\hat Q  \cap \{ Q_i \mid i =1,2,\dots \}$.  
 Then $Q = \hat Q*_{\hat F}Q'$, where $\hat F$ is
a locally finite subgroup of $Q$ which is a union of an increasing sequence of finite subgroups $F_i'$ where $F_i' \leq F_i$.
As in the example above,  the group $\hat Q$ must be finitely generated and one ended.
It is finitely generated because $Q$ is finitely generated,  and when we write a generating set for $Q$ as words given by the finite graph of groups 
decomposition just described then we will get a finite set of generators for $\hat Q$ by writing
each generator of $Q$ as a word in the elements of the vertex groups of the tree product and then taking those elements that are in $\hat Q$.
It will have to be one ended because if it split over a finite subgroup, then this decomposition will be induced by a similar decomposition of $Q$,
since a locally finite subgroup of $\hat Q$ must lie in a conjugate of one of the factors of the splitting.
The one ended subgraph of $W$ must, under a quasi-isometry, correspond to a one-ended subgraph of $Y$  which determines a special end.
The graph $Y$ has no such subgraph.   We have a contradiction.

We have proved that the locally finite graph $Y$ is quasi-isometric to a vertex transitive graph, but it is not quasi-isometric to
a Cayley graph.

\end{document}